\documentclass[10pt]{amsart}
\usepackage{amsfonts}
\usepackage{amssymb}
\usepackage{amsthm}
\usepackage[dvips]{graphicx}
\theoremstyle{plain}
\newtheorem{theorem}{Theorem}
\newtheorem{lemma}[theorem]{Lemma}
\newtheorem{definition}[theorem]{Definition}
\newtheorem{corollary}[theorem]{Corollary}

\renewcommand{\epsilon}{\varepsilon}

\newcommand{\tB}{\tilde{B}_n^0}

\begin{document}
\title{A Special Subgroup of the Surface Braid Group}
\author{D. Jeremy Copeland}
\maketitle

\section{Introduction}
Herein we prove that if $M$ is a compact oriented Riemann surface of genus $g$, and $M^{[n]}$ is the classifying space of $n$ distinct, unordered points on $M$, then the kernel of the map $\pi_1(M^{[n]})\to H_1(M)$ is generated by transpositions for sufficiently large $n$.  Specifically we prove the following theorem:

\begin{theorem}
If $M$ is a polyhedron of genus $g$ with $n$ faces such that no face is a neighbor of itself and no two faces share more than one edge, then $B_n^0=\ker(\pi_1(M^{[n]})\to H_1(M))$ is generated by the edge set.  Explicitly, the basepoint of $M^{[n]}$ may be chosen to be a marked point in the interior of each face, and each edge may be viewed as a transposition of the marked points on the faces it separates.
\end{theorem}

\section{Definitions}

This section contains most of the pertinent definitions.  Henceforth, $M$ will be a polyhedron of genus $g$ with $n$ faces such that no face is a neighbor of itself and no two faces share more than one edge.  If $S$ is any orientable surface, let $S^n$ be the subset of $S^{\times n}$ with all points distinct.  Fadell and Neuwirth have shown \cite{Fadell} that forgetting a coordinate gives a fibration $S^n\to S^{n-1}$ with fiber homeomorphic to $S\setminus\{s_1,\cdots,s_{n-1}\}$, where the $s_j$ are distinct points in $S$.  We will at some point use the homotopy extension and lifting principle of this map.

Let the set $\{f_j\}$ be the set of faces of $M$, $E$ its edges, and $V$ its vertices. Choose a basepoint $x_j$ for each face $f_j$.  Let $X=\{x_j\}$.  The edge set, $E$, of $M$ defines a graph $\Gamma$, and we have the dual graph $\check{\Gamma}$ with edge set $E$ and vertices $X$.  To each edge $e\in E$ viewed as an edge of $\check{\Gamma}$, we attach an orientation.  Thus any curve crossing $e\in\Gamma$ crosses with sign relative to $e\in\check\Gamma$.  We choose $X$ as the basepoint for all homotopy groups.

\begin{definition}
Let $x,y\in X$, and let $U\subset M$ open be such that $U\cap X=\{x,y\}$ and $U$ is simply connected.  
When we say \textbf{transposition} in $\pi_1(M^{[n]},X)$ of $x$ and $y$, we mean a generator of $\pi_1(U^{[2]},\{x,y\})$ extended by constants.  
\end{definition}

The groups in question are as follows:  
$\pi_1(M^{[n]})$ is the surface braid group.
$S^1=\mathbb{R}/ \mathbb{Z}$.
$A:\pi_1(M^{[n]})\to H_1(M)$ is the map which takes a curve $\gamma:S^1\to M^{[n]}$, $\gamma=\{\gamma_1,\ldots,\gamma_n\}$ to the union of the images $\gamma_j([0,1])$.  This gives a collection of closed curves without boundary.  
$B_n^0=\ker(A:\pi_1(M^{[n]})\to H_1(M))$.
If $e\in E$, $\partial e=x_j-x_k$, then a transposition of $x_j$ and $x_k$ across $e$ maps under $A$ to $e-e=0$.  Let $\tilde{B}_n^0$ be the subgroup of $B_n^0$ generated by the transpositions associated to $E$.  Clearly we intend to show that $\tilde{B}_n^0=B_n^0$.
For any surface $N$ we let $\Pi_1(N)$ be the fundamental groupoid of $N$.  

If $\gamma:([0,1],\{0,1\}) \to (M^{[n]},X)$ is a based map, we label $\gamma=\{\gamma_1,\ldots,\gamma_n\}$ by $\gamma_j(0)=x_j$.  We use $\gamma\mapsto[\gamma]$ to denote the map from based curves to $\pi_1(M^{[n]})$.

\begin{definition} Let $\gamma=\{\gamma_1,\ldots,\gamma_n\}$ be a curve in $M^{[n]}$ with $\gamma(0)=\gamma(1)= X$, which avoids the set $V$. Assume that $\gamma$ intersects $E$ for only finitely many $t$, and at each intersection, $\gamma_j(t)\in E$, $\gamma_j(t-\epsilon)$ and $\gamma_j(t+\epsilon)$ are contained in different faces for small $\epsilon$.  We call such a $\gamma$ \textbf{enumerable}.  
\end{definition}

On can show that for any $b\in B_n^0$, there is an enumerable $\gamma$ such that $[\gamma]=b$.  

\begin{definition} If $\gamma=\{\gamma_1,\ldots,\gamma_n\}$ is an enumerable curve in $M^{[n]}$, the \textbf{edge set of $\gamma_j$} is the ordered set of edges $E_{\gamma_j}=(e_{j1}^\pm\ldots e_{jp_j}^\pm)$ which are crossed by $\gamma_j$.  The sign is positive if $\gamma_j$ crosses $e_\alpha\in\Gamma$ in the same orientation as $e_\alpha\in\check\Gamma$.  Notice that the image of $\gamma_j$ in $\Pi_1(M\setminus V)=\Pi_1(\check{\Gamma})$ is the curve $e_{jp_1}^\pm\ldots e_{j1}^\pm$.  Define the \textbf{edge set of $\gamma$} as the collection of edge sets of its constituent curves: $E_\gamma=\{E_{\gamma_j}\}$.  We will write $E_j=E_{\gamma_j}$.
\end{definition}

The edge set $E_j$ is the signed sequence of edges encountered by $\gamma_j$.

\begin{definition}
If $\gamma=\{\gamma_1,\ldots,\gamma_n\}$ is an enumerable curve in $M^{[n]}$, and $E_{\gamma_j}=(e_{j1}^\pm\ldots e_{jp_j}^\pm)$ is the edge set of $\gamma_j$, then the \textbf{face set of $\gamma_j$} is the set $X_{\gamma_j}=(x_{j0}\ldots x_{jp_j})$ such that $e_{jk}^\pm$ is the edge from $x_{j,k-1}$ to $x_{jk}$.  Specifically, $x_{j0}=\gamma(0)$,  $x_{jp_j}=\gamma(1)$.
\footnote{Notice that the edge and face sets are ordered left to right, but multiplication of curves is right to left.}
\end{definition}

The face set $X_j$ is the sequence of faces encountered by $\gamma_j$.

\begin{definition}
We say that an edge set $X_j=(x_{j0}\ldots x_{jp_j})$ is a \textbf{palindrome} if $p_j$ is even and $x_{jk}=x_{j(p_j-k)}$ for all $k$.  We will say that an edge set is a \textbf{palindrome} if its corresponding face set is a palindrome.
\end{definition}

\begin{definition} 
We say that an enumerable curve $\gamma$ is \textbf{balanced} if for every edge $e_\alpha$, the multiplicity of $e_\alpha^+$ in $E_\gamma$ is equal to the multiplicity of $e_\alpha^-$.
\end{definition}

We are interested in balanced curves for the following reason:

\begin{theorem}
For any element $b\in B_n^0$, there is some balanced curve $\gamma$ such that $[\gamma]=b$.  
\end{theorem}

\begin{proof}
Choose any enumerable $\gamma'$ with $[\gamma']=b$.  $\gamma'$ need not be balanced, but since $b\in B_n^0$, the edge set of $\gamma'$ is exact (a boundary of the homology complex of $\check\Gamma$: $0\to \mathbb{Z}V\to  \mathbb{Z}E\to  \mathbb{Z}X\to 0$).  However, for any vertex $v$ (face of $\check{\Gamma}$), there are clearly curves $\gamma''$ such that $[\gamma'']=1$, and $E_{\gamma''}=\partial v$.  For example, if $\gamma_j''\equiv x_j$ for $j\geq 2$, and $\gamma_1''$ travels from $x_1$ to $v$, once around $v$, and returns along the original path to $x_1$, then $A\gamma''=\pm\partial v$.  Thus we may choose $\gamma''$ such that $[\gamma'']=1$ and $E_{\gamma''}=E_{\gamma'}$, so that $\gamma=(\gamma'')^{-1}\gamma'$ is the desired curve.

\end{proof}

\section{One particle motion}

In this section, we endeavor to prove that all one particle motion in $B_n^0$ lies in $\tB$.
First we will prove the following theorem.

\begin{theorem}\label{thm.ind}
Let $\gamma=\{\gamma_1,\ldots,\gamma_n\}:S^1\to (M^{[n]},X)$ such that $\gamma_j$ is constant for $j>1$ and $X_1$ is a palindrome, $(x_{j_1}x_{j_2}\ldots x_{j_{m-1}}x_{j_m}x_{j_{m-1}}\ldots x_{j_2}x_{j_1})$.  Then $[\gamma]\in \tilde{B}_n^0$.    
\end{theorem}

We prove this theorem by induction on $m$.  We call $m$ the \textbf{height} of $\gamma$.

\begin{lemma}\label{lem.lem1}
If $\gamma$ is as in Theorem \ref{thm.ind} and for some $k<m$, $x_{j_{k-1}}=x_{j_{k+1}}$, then there exist $\gamma'$, $\gamma''$, and $\gamma'''$ each with shorter palindromic face sets such that $[\gamma]=[\gamma''][\gamma'''][\gamma']$.  
\end{lemma}

\begin{proof}
We give two constructions for this proof.
Assume that for some $k<m$, $x_{j_{k-1}}=x_{j_{k+1}}$.  

Construct $\gamma'$ which follows $\gamma$ until reaching $x_{j_k}$, wrapping an appropriate number of times around $x_{j_m}$, then retracing its path to $x_1$.  Likewise construct $(\gamma'')^{-1}$ following $\gamma^{-1}$ in the same way (recall that $E_1$ is a palindrome, though $\gamma$ and $\gamma^{-1}$ need not behave in the same way at each vertex).  Now, there is some curve in the homotopy class $[\gamma'']^{-1}[\gamma][\gamma']^{-1}$ with only its first coordinate non-constant and such that its face set is $X_1$ with  $x_{j_{k-1}}$ and $x_{j_{k}}$ removed.  Let this be $\gamma'''$.

An alternative way to see this is as follows.  If $x_p$ and $x_q$ are neighbors, then at any time $\gamma_1$ travels through $x_p$, a homotopic curve may be obtained by following the same path from $x_1$ to $x_p$, then travelling to  $x_q$, returning to $x_p$, and completing $\gamma$.  The effect this has on the face set is $(\ldots x_p \ldots)\mapsto (\ldots x_p x_q x_p \ldots)$.  Beginning at $p=j_{k+1}$ we see that $[\gamma]$ may be realized by a curve with
\[
E_1=(x_{j_1}\ldots x_{j_{k-1}}x_{j_k} x_{j_{k-1}}\ldots  x_{j_2} x_{j_1} x_{j_2}\ldots x_{j_{k-2}} x_{j_{k+1}}\ldots x_{j_m}\ldots x_{j_1})
\]
Since this curve returns to $x_1$, one may factor it into two curves, the first being $\gamma'$.  Factoring $\gamma''$ from the end of $\gamma$ in the same way, we arrive at the curve $\gamma'''$.
\end{proof}

\begin{proof}[Proof of Theorem \ref{thm.ind}]
For simplicity we will assume that $j_1=1$, $j_2=2$, and that $e_\alpha$ connects $x_1$ to $x_2$.  Let $[\sigma_\alpha]$ be either transposition of $x_1$ and $x_2$.  We will write $\sigma_\alpha$ for some specific curve in $M^{[n]}$ fixing $x_j$, $j>2$ and $[\sigma_\alpha]$ for the element in $B^0_n$.
$\sigma_\alpha^2$ is homotopy equivalent to the curve in which $\gamma_1$ travels from $x_2$, once around $x_1$ and back to $x_2$, $\gamma_2\equiv x_1$, and $\gamma_j\equiv x_j$,  $j>2$.

$\sigma_\alpha\gamma\sigma_\alpha^{-1}$ is a curve in which $\gamma_1$ travels from $x_1$ across $e_\alpha$ to $x_2$, waits, then returns along $e_\alpha$ to $x_1$.   Thus $\gamma_1$ follows a homotopically trivial path, while $\gamma_j$ remain constant for $j>2$.  Using the Fadell-Neuwirth fibration $\check{\Gamma}^n\to\check{\Gamma}^{n-1}$, forgetting the coordinate of $x_2$, we may lift this homotopy to get a curve $\gamma'$ with $[\gamma']=[\sigma_\alpha][\gamma][\sigma_\alpha]^{-1}$, such that $\gamma_j$ is constant for $j>1$ and $\gamma_1(0)=x_2$.

We will construct $\gamma'$ explicitly.  Consider the map, $\gamma_1:[0,1]\to\check\Gamma$.  We have realized this curve as a curve which travels from vertex to vertex, turning some number of times around each vertex.  Lifting the homotopy has the effect of dragging the vertex $x_2$ to $x_1$.  Thus every time $\gamma_1$ reaches $x_2$ aside from the first and last, $\gamma'_1$ reaches $x_2$ travels to $x_1$ and around some number of times, then returns to $x_2$.  
\bigskip

\ \ \ \ \ \ \
\includegraphics[width=.3\textwidth]{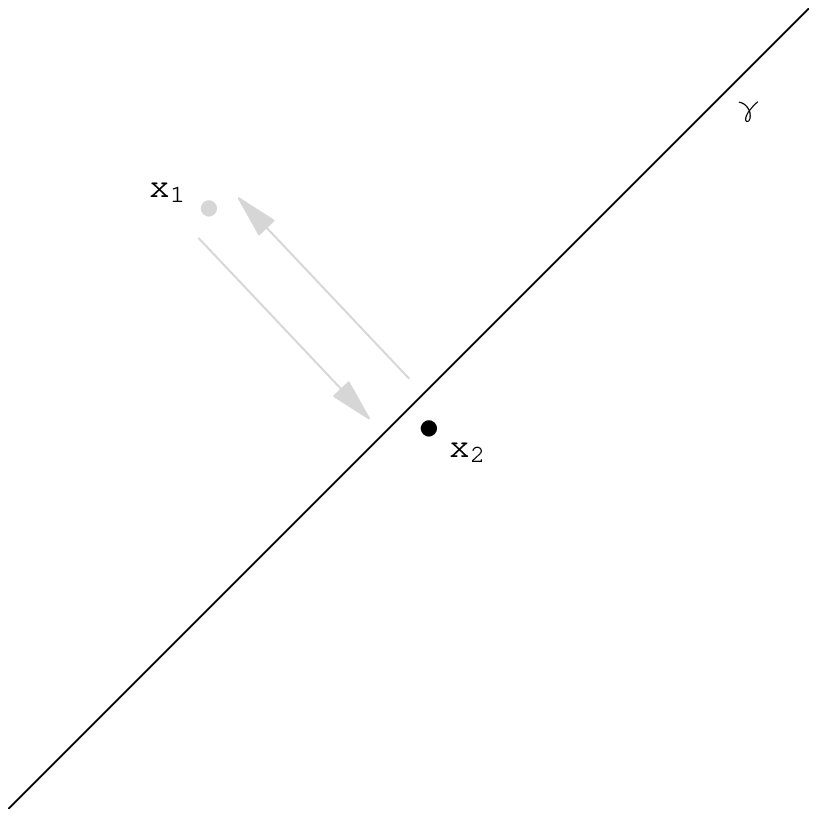}
\ \ \ \ \ \ \ \ \ \ \ \ \ \ \ \ \ \ \ \ 
\includegraphics[width=.3\textwidth]{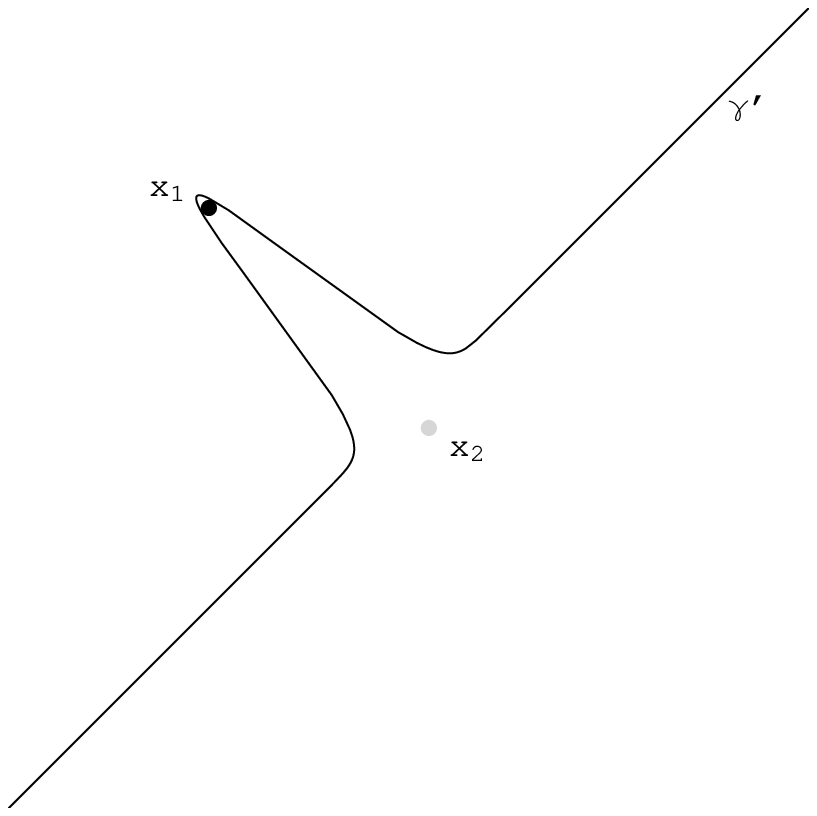}

\noindent 
Figure 1. $\gamma$ becomes $\gamma'$ under conjugation. 

Notice that there is some $\epsilon>0$ such that $\gamma_1|_{[\epsilon,1-\epsilon]}$ is a curve in $\hat\Gamma\setminus\{x_j|j\neq1\}$ while $\gamma'_1$ is a curve in  $\hat\Gamma\setminus\{x_j|j\neq2\}$.  Thus we have drawn $x_1$ or $x_2$ in light gray in our figures as points that $\gamma_1$ or $\gamma'_1$ may cross without affecting its homotopy class.

By construction, $E_{\gamma_1'}=(e_\alpha^+e_\alpha^-\cdots e_\alpha^+e_\alpha^-)$.  Therefore $\gamma'_1$ begins and ends by circling $x_1$ some number of times and returning to $x_2$.  This means that there is some $\gamma''$ and integers $a$ and $b$ such that 
\[
[\gamma'']=
[\sigma_\alpha]^{2b}[\gamma'][\sigma_\alpha]^{2a}=
[\sigma_\alpha]^{2b+1}[\gamma][\sigma_\alpha]^{2a-1},
\]
and that $E_{\gamma_1''}$ is the same as  $E_{\gamma_1'}$ without the beginning and ending copies of $e_\alpha^+e_\alpha^-$.  Left and right multiplication of $\gamma$ by odd powers of $\sigma$ ($L_{\sigma^{2b+1}}\circ R_{\sigma^{2a-1}}$) modifies the face set.  It acts by taking each $x_2$ to $x_2x_1x_2$ (as in Figure 1).  This is clear because this conjugation is the homotopy lifting map $\pi_1(\check{\Gamma}\setminus\{ \hat{x_1},x_2,\ldots,x_n\})\to \pi_1(\check{\Gamma}\setminus\{ x_1,\hat{x_2},\ldots,x_n\})$ dragging $x_2$ along $e_\alpha^-$ to $x_1$.   

We have constructed a new curve, $\gamma''$ such that $X_{\gamma_1''}$ is similar to $X_{\gamma_1}$ except that the first and last faces, which were both $x_1$, have been removed, so that the first and last faces are now $x_2$ and all other instances in $E_{\gamma_1}$ of $x_2$ have been replaced by the string $x_2x_1x_2$.  Also, $\gamma''_j$ are constant for $j>1$.
\[
(x_1x_2\ldots x_2\ldots x_2x_1)\mapsto (x_2\ldots x_2x_1x_2\ldots x_2).
\] 
Notice that this new edge set is still a palindrome.

Now we may apply Lemma \ref{lem.lem1} to factor $\gamma''$ into a product of curves.   Continue to use $x_{j_k}$ to denote the faces associated to $\gamma$.  For every $k>2$ such that $j_k=2$, $X_{\gamma''_1}$ contains a copy of the string $x_2x_1x_2$.  When $k\neq m$, we apply Lemma \ref{lem.lem1} to remove this string from the left and right.  Working from smallest $k$ to largest, we remove palindromes of height $k$ from left and right, replacing instances of $x_2x_1x_2$ in the middle palindrome with $x_2$.

We thus write $\gamma''$ as a product of $2K+1$ curves, where $K$ is the number of $2<k<m$ such that $j_k=2$.  For each such $k$, there are two curves of height $k$, and one final curve of height $m$ or $m-1$.  This final curve is the middle curve, and has the same face set as $\gamma$, except removal of the first and last faces, and if $j_m=2$, the middle face is replaced: $x_2\mapsto x_2x_1x_2$.  Thus if $j_m\neq 2$, then the final curve has height $m-1$ and we are done by induction on $m$.  

Assume that $j_m=2$. $[\sigma_\alpha]^{2b+1}[\gamma][\sigma_\alpha]^{2a-1}$ factors into a set of curves, but the central curve has height $m$.  That is, after applying Lemma \ref{lem.lem1}, we factor out curves of smaller height, and the remaining curve has face set 
\[
X_1=(x_{j_2}x_{j_3}\ldots x_{j_m}x_{1}x_{j_m}\ldots x_{j_3}x_{j_2}).
\]
  If $j_3\neq 1$ then we are done as before by applying the algorithm to get a palindrome of height $m-1$.  In the case that $j_3=1$, apply the argument again to get a curve with edge set 
\[
E_1=(x_{j_3}\ldots x_{j_m}x_{j_1}x_{j_2}x_{j_1}x_{j_m}\ldots x_{j_3})=(x_1\ldots x_2x_{1}x_2x_1x_2\ldots x_1).
\]
This again has height $m$.  However, notice that part of this curve travels through $x_2x_1x_2$.  This part of the curve is an element of $\Pi_1(f_1\cup f_2\setminus \{x_2\})$ based in $f_2$.  This is contractible to $f_2\setminus \{x_2\}$, so we may reduce the face set by $x_2x_1x_2\mapsto x_2$.  That is, a homotopy of $\gamma_1$ may pass freely through $x_1$, since $\gamma_1$ is a curve in $\check\Gamma\setminus\{x_j|j>1\}$.  When doing this, we create a new curve of smaller height which is still a palindrome.  Thus the theorem is proved.
\end{proof}

\begin{theorem}\label{thm.main}
If $\gamma=\{\gamma_1,\ldots,\gamma_n\}:S^1\to (M^{[n]},X)$ such that $[\gamma]\in B_n^0$, and $\gamma_j\equiv x_j$ for $j>1$, then $\gamma\in \tilde{B}_n^0$.    
\end{theorem}

\begin{proof}
The face set $X_1$ is not necessarily a palindrome (or, in fact, finite).  We need only factor some representative of $[\gamma]$ into curves for which the face sets are palindromes.   Assume, without loss of generality, that $\gamma$ has finite face set.  If $\phi:\tilde{M}\to M$ is the universal cover of $M$, let $\hat{M}$ be a finite sub-polyhedron with boundary in $\tilde{M}$ which contains $\phi^*\gamma$.  In fact, $\hat\pi=\pi_1(\hat{M}\setminus\phi^{-1}(\{x_2,\ldots,x_n\}))$ maps to $B^0_n$ and $[\gamma]$ is in its image.  However, it is well-known that $\hat\pi$ is generated by curves which travel from $x_1$ to $x_j$, pass around $x_j$, and return to $x_1$.  Such a curve is a palindrome.
\end{proof}

\section{Many particle motion}

Ultimately, we intend to show that $B_n^0=\tilde{B}_n^0$.  Since we may represent elements of  $B_n^0$ by balanced curves, we will study equivalences of balanced curves modulo $\tilde{B}_n^0$.  The first such equivalence says that if any curve begins by crossing $e_\alpha^+$ then $e_\alpha^-$, then it is equivalent to a curve which had not moved at all.

\begin{lemma}\label{lem.1}
If $\gamma$ is a balanced curve in $B_n^0$ with edge set $\{E_j\}$, and \\ $E_1=(e_\alpha^\pm e_\alpha^\mp e_{13}e_{14}\ldots e_{1p_1})$, then there is some element $b\in \tilde{B}_n^0$ and balanced $\gamma'$ such that $[\gamma]=[\gamma']b$ and the edge set of $\gamma'$ is $\{E_1',E_2,\ldots,E_n\}$ with $E_1'=(e_{13}e_{14}\ldots e_{1p_1})$.
\end{lemma}

\begin{proof}
Assume that $\gamma_1(t_0)=x_1$ and that $\gamma_1|_{[0,t_0]}$ passes through only the edges $e_\alpha^+$ and $e_\alpha^-$.  Assume further that no other curve passes through $x_1$.  We may vary $\gamma$ slightly to insure these assumptions.  Construct a new curve $\gamma'$ by $\gamma'_j=\gamma_j$ for $j>1$ and 
\[
\gamma'_1(t)=\left\{
\begin{matrix}
x_1&t<t_0\\
\gamma(t)&t>t_0.
\end{matrix}\right.
\]
Notice that $\gamma'$ has the appropriate edge set.  Consider the curve $\delta=(\gamma')^{-1}\gamma$ under the Fadell-Neuwirth fibration forgetting the first coordinate.
\footnote{As $\gamma_2,\ldots\gamma_n$ are constant, there is no problem lifting from $\check\Gamma^{[n]}$ to $\check\Gamma^n$.}
Clearly the image of this curve is trivial.  Therefore there is some element of $[\gamma']^{-1}[\gamma]$ such that all but the first coordinate is constant.  By Theorem \ref{thm.main} we are done.
\end{proof}

The second equivalence says that within an equivalence class, we may move the first letter of any word in the edge set to the beginning of some other word in the edge set.  e.g. 
\[
\{(play),(spies),(mexico),\ldots\}\equiv\{(splay),(pies),(mexico),\ldots\}
\]

\begin{lemma}\label{lem.edge}
If $\gamma$ is a balanced curve in $B_n^0$ with edge set $\{(e_{jk})\}$ and  $e_{11}=e_\alpha$ is an edge between $x_1$ and $x_2$, then there is some element $b\in \tilde{B}_n^0$ and balanced $\gamma'$ such that $[\gamma]=[\gamma']b$ and the edge set of $\gamma'$ is $\{(e_{12}\ldots e_{1p_1}),(e_\alpha e_{21}e_{22}\ldots e_{2p_2}\ldots),E_3,\ldots,E_n\}$.
\end{lemma}

\begin{proof}
Prepend to $\gamma$ the curve $\sigma_\alpha^{-1}$ so that $[\gamma]=[\gamma\sigma_\alpha^{-1}][\sigma_\alpha]$.  $\gamma\sigma_\alpha^{-1}$ now has edge set $\{(e_\alpha^{-1}e_\alpha e_{12}\ldots e_{1p_1}),(e_\alpha e_{21}e_{22}\ldots e_{2p_2}\ldots),E_3,\ldots,E_n\}$, where $E_j$ is as in $\gamma$ for $j>2$, since only if $\gamma_j(0)=x_1,x_2$ will the edge set change.  Finally one applies Lemma \ref{lem.1} to $[\gamma\sigma_\alpha]$.
\end{proof}

This obviously gives an equivalence of face sets:

\begin{corollary}\label{cor.mvface}
In the situation of Lemma \ref{lem.edge}, if the face set of $\gamma$ is $\{(xy\ldots)(y\ldots)\ldots\}$, then the face set of $\gamma'$ is identical except for the switch $\{(y\ldots)(xy\ldots)\ldots\}$.
\end{corollary}

Note that $x$ must appear at the beginning of exactly one word, call it the $k ^{\text{th}}$, and likewise with $y$.  This theorem allows us to move $x$ from the beginning of the $k ^{\text{th}}$ to the beginning of the word starting with $y$.  We will refer to this procedure as ``removing the first face from the $k^{\text{th}}$ word.''  

We are now ready to prove the main theorem.

\begin{theorem}
If $M$ is a polyhedron of genus $g$ with $n$ faces such that no side is a neighbor of itself and no two sides share more than one edge, then $B_n^0=\ker(\pi_1(M^{[n]})\to H_1(M))$ is generated by the edge set.  Specifically, the basepoint of $M^{[n]}$ may be chosen to be a marked point in the interior of each face, and each edge may be viewed as a transposition of the marked points on the faces it separates.
\end{theorem}

\begin{proof}
We show that any balanced curve $\gamma$ is equivalent modulo $\tilde{B}_n^0$ to a balanced curve $\gamma'$ with a smaller face set.  It is clear that the only curve with face set of size zero is the identity.  

We will do this by repeatedly applying Corollary \ref{cor.mvface} to remove the first face from various face sets in $\gamma$ until we reach the case $X_1=(x_2x_1x_2)$ or $X_1=(x_1x_jx_1x_2)$.  At this point, we apply Lemma \ref{lem.1} to reduce to $X_1=(x_2)$ or $X_1=(x_1x_2)$ while fixing all other $X_j$.  Thus the face set is smaller, and by induction, we are done.  In the proof of this theorem, we abandon the caveat that $X_j$ has first face $x_j$.

Let $X_1$ be some curve with a nontrivial face set.  First, we repeatedly remove the first face from the first curve so that $X_1$ has two elements.  Relabel the $x_j$ so that $X_1 =(x_1x_2)$.  Let $e_\alpha$ denote the edge connnecting $x_1$ to $x_2$.

Since $e_\alpha\in E_\gamma$ and $\gamma$ is balanced, some edge set contains the pair $x_2x_1$.  Call this edge set $X_2$.  Notice that the first instance of $x_1$ in $X_2$ need not occur as $\ldots x_2x_1\ldots$.  Next, we remove the first face $X_2$ until its first face is $x_1$.  At this point, $X_1=(x_jx_1x_2)$ and $X_2=(x_1x_k\ldots)$.  If $j=k$, then removing the first face of $X_2$ forces $X_1=(x_1x_jx_1x_2)$ and we are done.  Otherwise, removing the first faces of $X_2$ then $X_1$ places us in the situation $X_1=(x_1x_2)$ and $X_2=(x_k\ldots)$: 
\[
\begin{matrix}
X_1=(x_jx_1x_2)\\
\phantom{.}X_2=(x_1x_k\ldots)\\
X_?=(x_k\ldots)\ \ \ \\
\end{matrix}
~~\mapsto~~~
\begin{matrix}
X_1=(x_jx_1x_2)\\
X_2=(x_k\ldots)\phantom{...}\\
\phantom{}X_?=(x_1x_k\ldots)\\
\end{matrix}
~~\mapsto
\begin{matrix}
X_1=(x_1x_2)\phantom{.}\\
X_2=(x_k\ldots)\\
\phantom{aaa.}X_?=(x_jx_1x_k\ldots)\\
\end{matrix}
\]
In this process, we have reduced the number of times $x_1$ appears in $X_2$.  If we repeate this process exactly, eventually the first instance of $x_1$ in $X_2$ will be either $x_jx_1x_j$ or $x_2x_1$.  Then we are done. 
\end{proof}

\section{acknowledgements}
The idea of using the edges of a cellular decomposition as a generating set was suggested in a correspondence of R. Bezrukavnikov to V. Ginzburg, while Bezrukavnikov was a student at Tel-Aviv University.  
I'd like to thank my advisor, V. Ginzburg, for helpful suggestions and discussions.

\small{

}

\end{document}